\documentclass[12pt,a4paper,leqno]{amsart}

\usepackage{anysize}
\usepackage[utf8]{inputenc}
\usepackage{amsmath}
\usepackage{amssymb}
\usepackage{array}

\marginsize{3cm}{3cm}{3.5cm}{3.5cm}

\newtheorem{theorem}{Theorem}

\begin{document}

\title[Proof of the Ramanujan-type series $3A23$]{Proof of a rational Ramanujan-type series for $1/\pi$. The fastest one in level 3}

\author{Jesús Guillera} 

\address{University of Zaragoza, Department of mathematics, 50009 Zaragoza (Spain)}

\email{jguillera@gmail.com}

\dedicatory{To Bruce Berndt, in admiration of his inspirational work on Ramanujan's Notebooks}

\maketitle

\begin{abstract}
Using a modular equation of level $3$ and degree $23$ due to Chan and Liaw, we prove the fastest convergent rational Ramanujan-type series for $1/\pi$ of level $3$.  
\end{abstract}

\section{The method}
A prerequisite to understand this paper well is familiarity with the notation and the method developed in \cite{Gui-meth-rama}. This method is based on an original idea of Wan \cite{Wan}, whose paper was in turn influenced by some ideas of \cite{Chan-Wan-Zudilin}. 
Using the notation
\[
F_s(\alpha) = {}_2F_1\biggl(\begin{matrix} \frac1s, \, 1-\frac1s \\ 1 \end{matrix}\biggm| \alpha \biggr), \quad s \in \{ 6,4,3,2 \}, \qquad \ell=4\sin^2 \frac{\pi}{s},
\]
we proved in \cite{Gui-meth-rama} the following result:
\begin{theorem}
Let 
\[
F_s(\alpha) = m(\alpha,\beta) F_s(\beta), \quad A(\alpha, \beta)=0,
\]
be a transformation of modular origin and degree $1/d$, and let $\beta_0=1-\alpha_0$ be a solution of  $A(\alpha_0, \beta_0)=0$, and $m_0=m(\alpha_0, \beta_0)$. Then, if
\[
m_0 = \frac{1}{\sqrt{d}}, \quad \text {or} \quad m_0=\frac{\sqrt{4d-\ell}}{2d}+\frac{\sqrt{\ell}}{2d} \, i, 
\]
we have
\[
\sum_{n=0}^{\infty} \frac{\left(\frac12\right)_n\left(\frac1s \right)_n\left(1-\frac1s \right)_n}{(1)_n^3} (a+bn) \, z^n = \frac{1}{\pi},
\]
where 
\begin{equation}\label{b-positive-terms}
z=4 \alpha_0 \beta_0, \quad b=2 (1-2\alpha_0) \, \sqrt{\frac{d}{\ell}}, \qquad a=-2\alpha_0 \beta_0 \frac{m'_0}{\alpha'_0} \, \frac{d}{\sqrt{\ell}}.
\end{equation}
or
\begin{equation}\label{b-alternating-terms}
z=4 \alpha_0 \beta_0, \quad b=2 \, (1-2\alpha_0) \sqrt{\frac{d}{\ell}-\frac14}, \qquad a=-2\alpha_0 \beta_0 \frac{m'_0}{\alpha'_0} \, \frac{d}{\sqrt{\ell}},
\end{equation}
respectively. The $'$ means differentiation with respect to the variable that we choose as independent.
\end{theorem}

\section{The formula $(A, 3, 23)$}
In this paper we prove the Ramanujan-type formula
\begin{equation}\label{A-3-23}
\sum_{n=0}^{\infty}  \frac{\left(\frac12\right)_n\left(\frac13\right)_n\left(\frac23\right)_n}{(1)_n^3} (14151 n + 827) \frac{(-1)^n}{500^{2n}} = \frac{1500 \sqrt{3}}{\pi}.
\end{equation}
It is the formula $(A, 3, 23)$ in the notation of \cite{Gui-meth-rama}, we will also refer to it as $3A23$, and is the fastest convergent rational Ramanujan-type series for $1/\pi$ of level $\ell=3$. It was discovered by Chan, Liaw and Tan \cite[eq. 1.19]{Chan-Liaw-Tan}. As $z_0=4\alpha_0\beta_0$ where $\beta_0=1-\alpha_0$, we get
\[
\alpha_0=\frac12 - \frac{53\sqrt{89}}{1000}, \quad \beta_0=1-\alpha_0 = \frac12 + \frac{53\sqrt{89}}{1000},
\]
Then, with a numerical approximation of $20$ digits, we have
\[
m_0=\frac{F_3(\alpha_0)}{F_3(\beta_0)} \approx 
0.20508654634905660459+0.037653278425410375946 \, i,
\]
which we identify as
\[
m_0 = \frac{\sqrt{89}}{46}+\frac{\sqrt{3}}{46} \, i.
\]
Observe that $m_0$ is of the form
\[
m_0 = \frac{\sqrt{4d-\ell}}{2d}+\frac{\sqrt{\ell}}{2d}\, i, \qquad |m_0|=\frac{1}{\sqrt{d}},
\]
with $d=23$. Hence, for proving (\ref{A-3-23}) with our method, we need a transformation of degree $1/d$ with $d=23$ for the level $\ell=3$, and with such a transformation we can prove it rigorously. This is done in next section using a modular equation in which the algebraic relation $A(\alpha, \beta)=0$ is written using two auxiliary variables $u$ and $v$, in the following way:  
\begin{equation}\label{mod-eq-type}
u^k=\alpha \beta, \quad v^k=(1-\alpha)(1-\beta), \qquad P(u, v)=0,
\end{equation}
where $P(u, v)$ is a polynomial in $u$ and $v$, and $k$ is a positive integer.

\section{The proof of $(A, 3, 23)$}
We use a modular equation of level $3$ and degree $23$ due to Heng Huat Chan and Wen-Chin Liaw \cite[Corollary 3.7]{Chan-Liaw}, in which we have replaced $u$ and $v$ with $u^2$ and $v^2$ respectively.
It is
\begin{equation}\label{u-v-al-be}
u^{12}=\alpha \beta, \quad v^{12}=(1-\alpha)(1-\beta), \quad P(u, v)=0,
\end{equation}
where
\begin{multline}\label{mod-eq-deg-23}
P(u, v) \! = \! (u^8+v^8)-12\sqrt{3} (u^7 v+u v^7)-87(u^6 v^2 + u^2 v^6)-84 \sqrt{3} (u^5 v^3+u^3 v^5) \\ -160(u^4 v^4)-2(u^4 + v^4)-15 \sqrt{3}(u^3 v + u v^3) - 48 (u^2 v^2) +1.
\end{multline}
If we let $\beta=1-\alpha$, then we see that $u^{12}=v^{12}$. If we choose 
\begin{equation}\label{uv}
v=\left( \frac{\sqrt{3}}{2}-\frac12 \, i \right) u,
\end{equation}
and replace it in $P(u,v)=0$, we have the equation
\begin{equation}
250(1+\sqrt{3} \, i)u^8 -\frac{95}{2}(1-\sqrt{3} \, i)u^4+1 = 0,
\end{equation}
which factors as
\[
\frac{1+\sqrt{3}\, i}{16}(10 u^4 + 1+\sqrt{3}\, i)(20u^2-\sqrt{3}-i)(20u^2 + \sqrt{3} + i )=0.
\]
One solution is
\[
u_0=\frac{\sqrt{15}+\sqrt{5}}{20}+\frac{\sqrt{15}-\sqrt{5}}{20} \, i, \quad v_0=\frac{\sqrt{15}+\sqrt{5}}{20}-\frac{\sqrt{15}-\sqrt{5}}{20} \, i.
\]
Then, from
\[
\alpha_0(1-\alpha_0) = u_0^{12}= -10^{-6}, \quad \beta_0=1-\alpha_0,
\]
we get
\[
\alpha_0=\frac12-\frac{53}{1000} \sqrt{89}, \quad \beta_0=\frac12+\frac{53}{1000} \sqrt{89}.
\]
Differentiating $P(u, v)=0$ with respect to $u$ at $u=u_0$, we find
\[
v'_0=\frac{-294573 \sqrt{3}+82573 \, i}{516854}.
\]
Then, differentiating $P(u, v)=0$ twice with respect to $u$ at $u=u_0$, we get
\[
v''_0=\frac{8674041040500000 \sqrt{5}(1-i)+3034180783431000 \sqrt{15}(1+i)}{17258921684500483}.
\]
From (\ref{u-v-al-be}), we see that
\begin{equation}\label{al-be}
u^{12}=\alpha \beta, \quad u^{12}-v^{12}+1 = \alpha + \beta.
\end{equation}
Differentiating (\ref{al-be}) with respect to $u$ at $u=u_0$, we obtain
$\alpha'_0$ and $\beta'_0$. Then, differentiating (\ref{al-be}) twice with respect to $u$ at $u=u_0$ we obtain $\alpha''_0$
and $\beta''_0$. As the multiplier is given by 
\begin{equation}\label{multiplier}
m^2=\frac{1}{d} \frac{\beta(1-\beta)}{\alpha(1-\alpha)} \frac{\alpha'}{\beta'},
\end{equation}
see \cite{Gui-meth-rama}, replacing the already known values at $u=u_0$, we get
\begin{equation}\label{m0}
m_0=\sqrt{\frac{1}{23} \frac{\alpha'_0}{\beta'_0}}=\frac{\sqrt{89}}{46}+\frac{\sqrt{3}}{46} i = \frac{\sqrt{4d-\ell}}{2d}+\frac{\sqrt{\ell}}{2d} \, i.
\end{equation}
Taking logarithms in (\ref{multiplier}), differentiating with respect to $u$, and dividing by $\alpha'$, we get
\[
\frac{m'}{\alpha'} = \frac{m}{2 \alpha'} \left( \frac{\beta'}{\beta}-\frac{\beta'}{1-\beta} -\frac{\alpha'}{\alpha} + \frac{\alpha'}{1-\alpha} + \frac{\alpha''}{\alpha'} - \frac{\beta''}{\beta'}\right),
\]
and from it, we obtain
\[
\frac{m'_0}{\alpha'_0}=\frac{827000}{69}.
\]
Finally, using the formulas
\begin{equation}
z=4\alpha_0 \beta_0, \qquad b=2 \, (1-2\alpha_0) \sqrt{\frac{d}{\ell}-\frac14}, \qquad a=-2\alpha_0 \beta_0 \frac{m'_0}{\alpha'_0} \, \frac{d}{\sqrt{\ell}},
\end{equation}
we obtain
\[
z=\frac{-1}{500^2}, \qquad b=\frac{4717}{1500} \sqrt{3}, \qquad a=\frac{827}{4500} \sqrt{3},
\]
and we are done. 

\section{Ramanujan-type series for $1/\pi$ and modular equations}
In a completely similar way we can prove other Ramanujan-type series for $1/\pi$ using modular equations \cite{Be-vol-3, Be-vol-5, Be-Bh-Ga, Chan-Liaw}, written in the form (\ref{mod-eq-type}). For example, for proving $(A, 2, 7)$, we can use the modular equation
\[
u^8=\alpha \beta, \quad v^8=(1-\alpha)(1-\beta), \quad P(u,v)=0, 
\]
where
\[
P(u,v)=(u^4+v^4)+8\sqrt{2}(u^3 v + u v^3) + 20 u^2 v^2-1,
\]
see \cite[Theorem 10.3]{Be-Bh-Ga, Be-vol-5}. For proving $(A, 3, 11)$, we can use
\[
u^{12}=\alpha \beta, \quad v^{12}=(1-\alpha)(1-\beta), \quad P(u,v)=0,
\]
where
\[
P(u, v)=(u^4+v^4) + 3\sqrt{3} (u^3 v + u v^3) + 6 u^2 v^2 -1.
\]
see \cite[Theorem 7.8]{Be-Bh-Ga,Be-vol-5}, and we can get the proofs of the formulas $(A, 3, 5)$ and $(P, 3, 5)$,  with the modular equation 
\[
u^6=\alpha \beta, \quad v^6=(1-\alpha)(1-\beta), \quad u^2 + v^2 +3u v -1 = 0.
\]
see  \cite[Theorem 7.6]{Be-Bh-Ga, Be-vol-5}. Note that all the above examples correspond to rational series. Of course one can use the same method to prove irrational instances.

\section{Conclusion}
As we pointed out in  \cite{Gui-meth-rama}, the main aspect of our method is that the modular
equations used to prove the Ramanujan-type alternating series for $1/\pi$ have a much lower degree than those in other methods. It would be interesting to analyze in depth this phenomenon. The reader is invited to compare the values of ${\emph d}$ in the tables in \cite{Gui-meth-rama} with the respective values of ${\emph N}$ given in the tables in \cite{Chan-Cooper}. We hope that the method developed here will attract others to find, for example, new modular equations of levels $\ell \geq 5$ and apply them to prove non-hypergeometric series of Ramanujan-Sato type for $1/\pi$  \cite{Chan-Chan-Liu, Chan-Cooper}. Of course a generalization of the method would be needed for that.

\end{document}